\documentclass{amsart}
\usepackage{amscd} 
\usepackage{amsfonts} 
\usepackage{amssymb} 
\usepackage{latexsym} 
\input{diagrams}
\newcommand{\ncm}{\newcommand}

\def\C{\mathbb{C}\,}

\def\ME{\mathcal{E}}

\newtheorem{theorem}{Theorem}[section]
\newtheorem{prop}[theorem]{Proposition}
\newtheorem{lemma}[theorem]{Lemma}
\newtheorem{cor}[theorem]{Corollary}

\newtheorem{lem&def}[theorem]{Lemma \& Definition}
\newtheorem{definition}[theorem]{Definition}
\newtheorem{example}[theorem]{Example}

\ncm{\End}{\mbox{\rm End}\,}

\def\Hom{\mbox{\rm Hom}\,}

\def\Im{\mbox{\rm Im}\,}

\def\id{\mbox{\rm id}}

\def\into{\hookrightarrow}
\def\to{\rightarrow}

\def\o{\otimes}    
\def\bra{\langle}
\def\ket{\rangle}

\ncm{\rarr}[1]{\stackrel{#1}{\longrightarrow}}
\ncm{\larr}[1]{\stackrel{#1}{\longleftarrow}}
\def\cop{\Delta}

\def\eps{\varepsilon}

\def\du1{\hat 1}

\def\-1{_{(-1)}}
\def\0{_{(0)}}
\def\1{_{(1)}}
\def\2{_{(2)}}
\def\3{_{(3)}}

\def\|{\, | \,}

\def\du1{\hat 1}
\def\lact{\triangleright}
\def\ract{\triangleleft}
\def\op{^{\rm op}}

\begin{document}

\title{The Endomorphism Ring Theorem for Galois and D2 extensions}
\author{Lars Kadison}
\address{Matematiska Institutionen \\ G{\" o}teborg 
University \\ 
S-412 96 G{\" o}teborg, Sweden} 
\email{lkadison@c2i.net} 
\date{}
\thanks{The author thanks Korn\'el Szlach\'anyi, Gabriella B\"ohm and K.F.K.I.\  for
a hospitable visit to Budapest in 2005 and stimulating discussions.}
\date{} 

\begin{abstract} 
Let $S$ be the left bialgebroid $\End {}_BA_B$ over
the centralizer $R$ of  a right D2 algebra extension $A \| B$,
which is to say that its tensor-square is isomorphic as $A$-$B$-bimodules to
 a direct summand of a finite direct sum of $A$ with itself.  
We prove that its left endomorphism algebra 
 is a left $S$-Galois extension of
$A^{\rm op}$.  
As a corollary, 
endomorphism ring theorems for D2 and
Galois extensions are derived from the D2 characterization of Galois extension.  
We note the converse that a Frobenius extension satisfying a generator condition
is D2 if its endomorphism algebra extension is D2.  
\end{abstract} 
\maketitle

\section{Introduction}

\begin{figure}
$$\begin{diagram}
\ME \o_{\rho(A)} \ME &&  \rTo^{\beta}&& S \o_R \ME \\
\dTo^{\cong} && && \uTo_{\cong}    \\
\ME \o_A \ME & &  && \Hom ({}_BA \o_B A, {}_BA) \\
& \SE_{\cong}& & \NE_{\cong} & \\
&& \Hom({}_B\Hom(\ME_A, A_A), {}_BA) &&
\end{diagram}$$
\caption{Galois map $\beta$ factors through various isomorphisms
in Theorem~\ref{th-endogalext}.}
\end{figure}

Bialgebroids are generalized weak bialgebras over an arbitrary
 noncommutative base ring \cite{BW,KS,T}. As in the theory of bialgebras or 
weak bialgebras,
there are associated to a bialgebroid, module and comodule algebras, smash products and Galois
extensions \cite{Bo,fer, Mo, NV, Sz}.
These constructions all occur in the tower over a depth two (D2) extension $A \| B$,
where the noncommutative base ring is the centralizer $R$ of the ring extension: the extension is  D2 if its tensor-square is centrally projective w.r.t.\ $A$ as a natural
$B$-$A$-bimodule (left D2) and $A$-$B$-bimodule (right D2)  \cite{fer, KS}. An extra condition that the natural module $A_B$
is balanced or faithfully flat ensures that $A \| B$ is a right Galois extension w.r.t.\ the $R$-bialgebroid $T := (A \o_B A)^B$
\cite{fer}.  
In section~4 the main theorem of this paper shows that even without this condition the left endomorphism ring $\ME := \End {}_BA$ is a left  $S$-Galois
extension of the subalgebra of right multiplications $\rho(A)$.  The proof that the
Galois mapping is an isomorphism is  summarized in the commutative diagram in Figure~1. 

Endomorphism ring theorems date back at least to the 1950's with Kasch's theorem for Frobenius extensions,
Nakayama and Tzuzuku's theorem for $\beta$-Frobenius extensions, 
and similar theorems for QF-extensions by M\"uller \cite{Mu}:  Morita's  paper \cite{Mo} was a 
definitive paper for all such Frobenius extensions.  The general idea
of these theorems is that important properties such as the Frobenius property of an extension $B \subseteq A$ or a ring homomorphism $B \to A$
sometimes pass up to the extension $A \into \End A_B$ induced from left multiplication $\lambda(a)(x) := ax$.
Also the index of a Frobenius extension, or the Hattori-Stallings rank of the underlying projective
module, may be viewed as passing up unchanged to the endomorphism ring extension.  
Sometimes two properties are dual with respect to this shift of levels such as separability and splitness
of an extension, as shown by Sugano, so that
a split extension leads to separable endomorphism ring extension, and \textit{vice versa}
under suitable conditions. Finally possessing a certain Frobenius homomorphism compatible, or ``Markov,'' trace
passes up in an endomorphism ring theorem as well. 
 These four curiosities as endomorphism ring theorems had a stunning
application to topology in the 1980's when Jones iterated the endomorphism ring construction of a tower over special
split, separable Frobenius algebra extensions (type $II_1$ subfactors), finding a countable set of idempotents that satisfy braidlike relations,
which together with a Markov trace lead to the first knot and link polynomials since the classical Alexander polynomial (see \cite{NEFE} for an explanation of this point of view,
and \cite{BW,KT} for the related coring point of view).

Depth two extensions have their origins in finite depth $II_1$ subfactors.
An inclusion of finite-dimensional C$^*$-algebras $B \subseteq A$ can be recorded as a bicolored weighted
multigraph called a Bratteli inclusion diagram:  
the number of edges between a black dot representing an (isomorphism class of a)
simple module $V$ of $A$
and a white dot representing a simple module $W$ of $B$ is $\dim {\rm Hom}_B(V,W)$.  This can be recorded in an inclusion matrix
of non-negative integers, which corresponds to an induction-restriction table of
irreducible characters of a  subgroup pair $H < G$ if $A = \C G$ and $B=\C H$ are the group algebras.  

If we define the basic construction of a semisimple $\C$-algebra pair $B \subseteq A$ to be the endomorphism algebra
$\End A_B$ of intertwiners, we note that $B$ and $\mathcal{E} := \End A_B$ are Morita equivalent via bimodule ${}_{\mathcal{E}}A_B$ 
whence the inclusion diagram of the left multiplication inclusion $A \into \End A_B$ is reflection of the 
diagram of $B \subseteq A$.  Beginning with a subfactor $N \subseteq M$, we build the Jones tower using the basic construction
$$ N \subseteq M \subseteq M_1 \subseteq M_2 \subseteq \cdots $$
where $M_{i+1} = \End {M_i}_{M_{i-1}}$, then the derived tower of centralizers or relative commutants
are f.d.\ C$^*$-algebras,
$$ C_N(N) \subseteq C_M(N) \subseteq C_{M_1}(N) \subseteq C_{M_2}(N) \subseteq \cdots $$
The subfactor $N \subseteq M$ has finite depth if the inclusion diagrams of the derived tower stop growing and begin
reflecting at some point, depth $n$ where counting begins with $0$. For example, the subfactor has depth two
if $C_{M_2}(N)$ is isomorphic to the basic construction of $C_M(N) \subseteq C_{M_1}(N) $.
As shown in \cite{KS}, depth two is really a property of ring or algebra extensions if viewed as a property
of the tensor-square. One of the aims in such a generalization is to find an algebraic theorem corresponding
to the Nikshych-Vainerman Galois correspondence (see \cite{NV}) between intermediate subfactors of a depth two
finite index subfactor $M \| N$ and left coideal subalgebras of the weak Hopf algebra  $C_{M_1}(N)$.

In \cite[Theorem 6.1]{KS} we noted an endomorphism ring theorem for Frobenius D2 extensions, where 
the proof shows that a right D2 extension has a left D2 endomorphism ring extension
(although right and left D2 are equivalent for Frobenius extensions).  In this
paper, Corollary~\ref{cor-leftendo} of the main theorem shows that Frobenius may be removed
from this theorem.  We use a characterization of one-sided Galois extension for bialgebroids
in terms of the corresponding one-sided D2 property of extensions in \cite[Theorem 2.1]{fer} or 
Theorem~\ref{th-characterization} in this paper.
This characterization and our main theorem then yields endomorphism ring theorems for D2 extensions
and Galois extensions:
a right D2 (or Galois) extension has left D2 (or Galois) left endomorphism extension $\ME \| A\op$ (Corollary~\ref{cor-leftendo}). There is a formula for the Galois inverse mapping in eq.~(\ref{eq: beta inverse})
(not depending on an antipode as in e.g.\ \cite[below eq.\ (13)]{KN}).
Keeping track of opposite
algebras of D2 extensions and their associated bialgebroids, we extend this to a right endomorphism 
theorem in Corollary~\ref{cor-rightendo}. Asking for converses to the endomorphism ring theorem usually lead to more difficult questions:
we note however a simple proof that a D2 endomorphism algebra extension implies
that a one-sided split or generator Frobenius extension  is itself D2
(Theorem~\ref{th-converse-endo}), which answers \cite[Question 1]{KN}.  

Bialgebroids equipped with antipodes are Hopf algebroids, although there is a scientific discussion about
what definition to use \cite{BS, CM, Lu}.  Finding further examples would tend to throw more weight to 
one of these non-equivalent definitions.  
In Theorem~\ref{th-newex} we show that the bialgebroid $T$ of a D2 extension $A$ over a Kanzaki
separable algebra $B$ 
is a Hopf algebroid of the type in \cite[B\"{o}hm-Szlach\'anyi]{BS}.  The antipode is very
naturally given by a twist of $(A \o_B A)^B$ utilizing a symmetric separability element.  
These provide  further examples of non-dual Hopf algebroids, in contrast to the dual Hopf
algebroids $S$ and $T$ of a D2 Frobenius extension $A \| B$ \cite{BS}.  They are also Hopf algebroids
with no obvious counterparts among Lu's version of Hopf algebroid \cite{Lu}.

\section{Preliminaries on D2 extensions} 

The basic set-up throughout this paper is the following.
We work implicitly with associative unital algebras
over a commutative ground ring $K$ in the
category of $K$-linear maps where bimodules are $K$-symmetric.  An algebra
extension $A \| B$ is a unit-preserving algebra homomorphism $B \to A$,
proper if this is monomorphic. We often say that an extension $A \| B$ has
left or right property X (such as being finite projective, i.e.\
f.g.\ projective) if the natural module ${}_BA$ or
$A_B$ has the property X.
Let
$A \| B$ be an algebra extension with centralizer denoted by $R := C_A(B) = A^B$,
bimodule endomorphism algebra $S := \End {}_BA_B$ and $B$-central
tensor-square $T := (A \o_B A)^B$. 
  $T$ has an algebra structure induced
from $T \cong \End {}_A A \o_B A_A$ given by 
\begin{equation}
\label{eq: tee mult}
tt' = {t'}^1 t^1 \o t^2 {t'}^2, \ \ \ \ \ 1_T = 1 \o 1,
\end{equation}
where $t = t^1 \o t^2 \in T$ uses a Sweedler notation and suppresses
a possible summation over simple tensors. Note the homomorphism
$\sigma: R \to T$ given by $\sigma(r) = 1 \o r$ and
the anti-homomorphism $\tau: R \to T$ by $\tau(r) = r \o 1$:
the images of $\tau$ and $\sigma$ in all cases commute in $T$.  
We will work with the $R$-$R$-bimodule induced by $T_{\sigma,\, \tau}$, which is given
by 
\begin{equation}
\label{eq: r bimod on t}
 {}_RT_R:\ \ \ r \cdot t \cdot r' := t \sigma(r') \tau(r) = rt^1 \o t^2 r'
\end{equation}
for all $t \in T$, $r,r' \in R$. 

There is a category of algebra extensions whose objects are algebra homomorphisms $B \to A$ 
 and whose morphisms from $B \to A$ to $B' \to A'$ are commutative squares with algebra homomorphisms
$B \to B'$ and $A \to A'$, an isomorphism if the last two maps are bijective.  The notions
of left or right D2, balanced, finite projective or Galois extensions introduced below are isomorphism invariants in this category (e.g.,
left D2 involves applying eq.~\ref{eq: left d2 qb}).     

 Let $\lambda: A \into \End A_B$
denote the left multiplication operator and $\rho: A \into \End {}_BA$ the right multiplication operator.
Let $\ME$ denote $\End {}_BA$ and
note that $S \subseteq \mathcal{E}$, a subalgebra under the usual
composition of functions.  Note that $\lambda$ restricts to an algebra monomorphism $R \into S$ and $\rho$ restricts
to an anti-monomorphism $R \into S$; and the images commute in $S$ since $\lambda(r) \rho(r') = \rho(r') \lambda(r)$
for every $r \in R$. 
We work with the $R$-$R$-bimodule ${}_{\lambda,\, \rho}S$ given by
\begin{equation}
\label{eq: r bimod on s}
{}_RS_R: \ \ \ r\cdot \alpha \cdot r' := \lambda(r) \rho(r') \alpha = r \alpha(-) r'
\end{equation}
for all $\alpha \in S$, $r,r' \in R$. 

We have the notion of an arbitrary bimodule
being centrally projective with respect to a canonical bimodule: we say that
a bimodule ${}_AM_B$, where $A$ and $B$ are two arbitrary algebras,
is \textit{centrally projective w.r.t.\ a bimodule} ${}_AN_B$, if ${}_AM_B$ is
isomorphic to a direct summand of a finite direct sum of $N$
with itself; in symbols,
if ${}_AM_B \oplus * \cong \oplus^n {}_AN_B$.  This recovers the usual notion
of centrally projective $A$-$A$-bimodule $P$ if $N = A$, the natural bimodule.  

We recall the definition of a D2, or \textit{depth two}, algebra extension $A \| B$ as simply  that  its tensor-square
$A \o_B A$  be centrally projective w.r.t.\ the natural $B$-$A$-bimodule $A$ (left D2)
and the natural $A$-$B$-bimodule $A$ (right D2). 

 A very useful characterization of D2 extension is that an extension is left D2 if there exist finitely many paired elements (a left D2 quasibase) $\beta_i \in S$,
$t_i \in T$  such that
\begin{equation}
\label{eq: left d2 qb}
a \o_B a'  =  \sum_i t_i \beta_i(a)a'  
\end{equation}
and right D2 if there are finitely many paired elements (a right D2 quasibase) 
$\gamma_j \in S$, $u_j \in T$ such that 
\begin{equation}
\label{eq: right d2 qb}
 a \o_B a' =  \sum_j a \gamma_j(a')u_j
\end{equation}
for all $a, a' \in A$ \cite[3.7]{KS}: we fix this notation, 
especially for a right D2 extension.  
We note from these equations and the dual bases characterization of finite projective modules
that the modules ${}_RT$ and $S_R$ are finite projective in case $A \| B$ is right D2,
while $T_R$ and ${}_RS$ are finite projective in case $A \| B$ is left D2. 
(For in the case of $S_R$, for example, let $a = 1$ in eq.~(\ref{eq: right d2 qb}) and apply $\id \o \alpha$ for
any $\alpha \in S$. Note that $T \longrightarrow$ 
$\Hom (S_R, R_R)$
via $t \mapsto$ 
$ (\alpha \mapsto$ 
$ t^1 \alpha(t^2))$, in fact an isomorphism.)

Centrally projective algebra extensions, H-separable
extensions and f.g.\
Hopf-Galois extensions are some of the classes of examples of D2 extension. If $A$ and $B$ are the complex
group algebras corresponding to a subgroup $H < G$ of  a finite group, then $A \| B$ is D2 iff $H$ is
a normal subgroup in $G$ \cite{KK}. 
 In \cite{fer} we show  that finite weak Hopf-Galois extensions are left and
right D2. More generally, left Galois extensions for bialgebroids and their comodule algebras are left D2 \cite{fer}. 

A surprisingly large portion of the Galois theory for bialgebroids
accompanying a D2 extension may be built up from just a left
or right D2 extension \cite{fer}. Working with as few hypotheses
as necessary, we will therefore  work with right D2 extensions
in this paper.  
  The proposition below and Props.~\ref{prop-teecop} and~\ref{prop-cop} further below indicate how one 
theory for either left or right D2 extensions may cover
both cases via a duality.  Let $A\op$ denote the opposite algebra of an algebra $A$.
We briefly use the notation $S(A \| B)$ and $T(A \| B)$ to denote the $S$ and $T$ constructions
above with emphasis on their dependence on the algebra extension
$A \| B$. 

\begin{prop}
\label{prop-op}
The algebra extension $A \| B$ is right D2 if and only if its opposite algebra
extension $A\op \| B\op$ is left D2.
Moreover, there are canonical algebra isomorphisms $S(A \| B) \cong S(A\op \| B\op)$ and
$T(A \| B) \cong T(A\op \| B\op)$.
\end{prop}
\begin{proof}
Let $a \mapsto \overline{a}$ be the identity anti-isomorphism
$A \to A\op$.  Then ${}_AA_B \cong {}_{B\op}{A\op}_{A\op}$
are identified via this same mapping, as are
${}_A A\o_B A_B \cong {}_{B\op} A\op \o_{B\op} {A\op}_{A\op}$
by $a \o c \mapsto $ $\overline{a \o c} :=$ $ \overline{c} \o \overline{a}$.  We also note
an isomorphism $\End A_B \cong \End {}_{B\op}A\op$ by
$f \mapsto \overline{f}$ where $\overline{f}(\overline{a}) :=$ $\overline{f(a)}$; similarly the left endomorphism algebra of $A \| B$
is identified with the right endomorphism algebra of $A\op \| B\op$.
The intersection of left and right endomorphism algebras, 
 $S((A \| B)$ and $S(A\op \| B\op)$ are therefore isomorphic as algebras. 
The algebras $T(A \| B) \cong T(A\op \| B\op)$ via
$t \mapsto \overline{t}$ is readily checked
to be an isomorphism of algebras via eq.~(\ref{eq: tee mult}).
It follows that paired elements $\gamma_j \in S(A \| B)$,
$u_j \in T(A \| B)$ form a right D2 quasibase for $A \| B$ iff
$\overline{\gamma_j} \in S(A\op \| B\op)$, $\overline{u_j} \in T(A\op \| B\op)$ form a left D2 quasibase for $A\op \| B\op$.  
\end{proof}

We add two new characterizations of  a one-sided D2 extension $A \| B$
in terms of restriction, induction and coinduction between $A$-modules and
 $B$- or $R$-modules.

\begin{prop}
Let $A \| B$ be  an algebra extension and $C$ an algebra.  $A \| B$ is right D2 if
and only if either one of the conditions below is satisfied: 
\begin{enumerate}
\item ${}_RT$ is finite projective and for every $C$-$A$-bimodule ${}_CM_A$,
there is a \newline 
$T \o_K C$-$B$-bimodule isomorphism
$$ M \o_B A \stackrel{\cong}{\longrightarrow} M \o_R T $$
\item $S_R$ is finite projective and for every $C$-$A$-bimodule ${}_CM_A$,
there is a \newline 
 $T \o_K C$-$B$-bimodule isomorphism
$$ M \o_B A \stackrel{\cong}{\longrightarrow} \Hom (S_R, M_R). $$
\end{enumerate}
\end{prop}
\begin{proof}
Suppose $A \| B$ is right D2 with quasibases $u_j \in T$ and $\gamma_j \in S$.
To the mapping $M \o_R T \to M \o_B A$ defined by 
\begin{equation}
\label{eq: emtee}
m \o_R t \longmapsto mt^1 \o_B t^2
\end{equation}
we define the inverse mapping by
\begin{equation}
m \o_B a \longmapsto \sum_j m \gamma_j(a) \o_R u_j.
\end{equation}
The left $T$-module structures on $M \o_R T$ and $M \o_B A$ are given by
$t \cdot (m \o t') := m \o_R tt'$ and $t \cdot (m \o a) := mt^1 \o_B t^2 a$,
both of which commute with the left $C$-module structure.
It is clear from eq.~(\ref{eq: tee mult}) that the mapping in eq.~(\ref{eq: emtee})
is left $T$-linear.
  That the mappings are inverse follows from  $t^1\gamma_j(t^2) \in R$
for each $j$ and from applying eq.~(\ref{eq: right d2 qb}) combined with the well-definedness of
some obvious mappings.

Furthermore, to the mapping $M \o_B A \to$ $ \Hom (S_R, M_R)$ defined by
\begin{equation}
\label{eq: em a}
m \o_B a \longmapsto (\alpha \mapsto m \alpha(a))
\end{equation}
we define the inverse mapping by
\begin{equation}
g \longmapsto \sum_j g(\gamma_j)u_j^1 \o_B u_j^2,
\end{equation}
where we note from eq.~(\ref{eq: right d2 qb})
that $\sum_j \gamma_j \cdot (u^1_j \alpha(u^2_j)) = $ $\alpha$.
The left $T$-module structure on $ \Hom (S_R, M_R)$ is given by
$(t \cdot f)(\alpha) := f(t^1 \alpha(t^2 -))$ for each $\alpha \in S$,
where we make use of a right action of $T$ on $\End {}_BA$ defined below in eq.~(\ref{eq: ract by tee})
and restricted to $S$. It is clear that this module action commutes with the left $C$-module structure
on $\Hom (S_R, M_R)$ induced by ${}_CM$, and that the mapping in eq.~(\ref{eq: em a}) is left $T \o C$-linear.   

Conversely, let $C = A$ and $M = A$, the natural $A$-$A$-bimodule,
and $\Leftarrow$ follows from \cite[Theorem 2.1, (3),(4)]{KK} applied to the opposite
algebra extension. (No use is made of the left $T$-module structure in this argument.)
 This is also easy to see by arriving at the defining condition 
for right D2 ${}_A A \o_B A_B \oplus * \cong \oplus^n {}_A A_B$, using 
the similar finite projectivity condition  for  ${}_RT$ or $S_R$,
and the functor $A \o_R - $ from left $R$-modules to $A$-$B$-bimodules as
well as the functor $\Hom (- , A_R)$ from right $R$-modules to $A$-$B$-bimodules. 
\end{proof} 

A corresponding characterization of left D2 extensions is derivable  by the opposite
algebra extension technique introduced in Prop.~\ref{prop-op}. 


\section{Preliminaries on bialgebroids and Galois extensions}

In this section we review the theory of bialgebroids, Hopf algebroids, their comodule algebras and Galois extensions.
We provide an example  of the Hopf algebroid in \cite{BS}
drawn from separability and finite dimensional algebras.  

Recall that a left bialgebroid 
$R'$-bialgebroid $S'$ is first of all two algebras $R'$ and $S'$ with two commuting maps (``source'' and ``target'' maps)
$\tilde{s}, \tilde{t}: R' \to S'$, an algebra homomorphism and anti-homomorphism resp., commuting in the sense  that $\tilde{s}(r) \tilde{t}(r') = \tilde{t}(r') \tilde{s}(r)$ for all $r,r' \in R'$.  Secondly, it is an $R'$-coring 
$(S',\ \cop \!: S' \to \ $ $ S' \o_{R'} S', \eps: S' \to R')$  \cite{BW} w.r.t.\ the $R'$-$R'$-bimodule $S'$ given by 
($x \in S;\, r,r' \in R$)
\begin{equation}
\label{eq: uno}
r \cdot x \cdot r' =  \tilde{t}(r') \tilde{s}(r)x.
\end{equation}
 Third, it is  a generalized bialgebra (and generalized weak bialgebra) in the sense that we have the axioms 
 \begin{equation}
\label{eq: duo}
\cop(x) (\tilde{t}(r) \o 1_S) = \cop(x)(1_S \o \tilde{s}(r)),
\end{equation}
$\cop(xy) = \cop(x)\cop(y)$ (which makes sense thanks to the previous axiom), $\cop(1) = 1 \o 1$,
and  $\eps(1_{S'}) = 1_{R'}$ and
 \begin{equation}
\label{eq: tres}
\eps(xy) = \eps(x\tilde{s}(\eps(y)) ) = \eps(x\tilde{t}(\eps(y)) )
\end{equation}
  for all $x, y \in T', r,r' \in R'$. The $R'$-bialgebroid $S'$ is denoted by
$(S', R', \tilde{s}, \tilde{t}, \cop, \eps)$.  

A right bialgebroid is defined like a left bialgebroid with three of the axioms transposed;
viz., eqs.~(\ref{eq: uno}),~(\ref{eq: duo}) and (\ref{eq: tres}) \cite{KS}. The left
bialgebroid $(S', R', \tilde{s}, \tilde{t}, \cop, \eps)$ becomes a right bialgebroid
under the opposite multiplication and exchanging source and target maps 
$${S'}\op =({S'}\op, R', \tilde{t}, \tilde{s}, \cop, \eps)$$
and a left bialgebroid under the co-opposite comultiplication, opposite base algebra
and exchanging source and target maps 
$${S'}_{\rm cop} = (S', {R'}\op, \tilde{t}, \tilde{s}, {\cop}\op, \eps),$$
which has the opposite $R'$-$R'$-bimodule structure from the $S'$ one starts with \cite[section 2]{KS}.  A similar situation
holds for opposites and co-opposites of right bialgebroids. 

 In \cite[4.1]{KS} we establish that $S := \End {}_BA_B$ is a left bialgebroid over the centralizer $R$ with the $R$-$R$-bimodule
structure given by eq.~(\ref{eq: r bimod on s}).
The comultiplication $\cop_S: S \to S\o_R S$ is given by either of two formulas, using
either a left or right D2 quasibase, equal in case $A \| B$ is D2: 
\begin{eqnarray} 
\label{eq: left d2 comult}
\cop_S(\alpha) & := & \sum_i \alpha(-t_i^1)t_i^2 \o \beta_i \\
\label{eq: right d2 comult}
    \cop_S(\alpha)   & := & \sum_j \gamma_j \o u_j^1 \alpha(u_j^2 -). 
\end{eqnarray}
Since $S \o_R S \cong \Hom_{B-B}(A \o_B A, A)$ via $\alpha \o \beta \mapsto$ $ (a \o a' \mapsto \alpha(a)\beta(a'))$ \cite[3.11]{KS},
with inverse given by either $F \mapsto \sum_i F(- \o t^1_i) t^2_i \o_R \beta_i$ or
$\sum_j \gamma_j \o_R u^1_j F(u^2_j \o - )$, our formulas for the coproduct simplify greatly via eq.~(\ref{eq: right d2 qb})
to $\cop_S(\alpha)(a \o a') = \alpha(aa')$ (thus a generalized Lu bialgebroid).  The counit belonging to this coproduct
is $\eps_S: S \to R$ given by $\eps_S(\alpha) = \alpha(1_A)$.  It is then not hard to see
that $(S, \cop, \eps)$ is an $R$-coring, i.e., both maps are $R$-bimodule morphisms, that $\cop_S$ is coassociative and that $(\eps_S \o \id)\cop_S =$ $ \id_S = $ $(\id_S \o \eps_S) \cop_S$.  In \cite{fer} it is shown that a right D2 quasibase alone leads to
a left bialgebroid $S$ over $R$ with $S_R$ finite projective; similarly, a left D2 quasibase on $A \| B$
implies $S$ is a left finite projective left  $R$-bialgebroid.

The algebra $T$ defined above for any algebra extension is a right bialgebroid over the centralizer 
$R$.  The $R$-coring structure underlying the right $R$-bialgebroid $T$ is given by
an $R$-$R$-bimodule structure $r \cdot t \cdot r' = $ $ rt^1 \o t^2 r'$, a comultiplication equivalently  definable using either a right D2 quasibase or left D2 quasibase, 
\begin{eqnarray}
\label{eq: tee}
\cop_T(t) & := & \sum_j (t^1 \o_B \gamma_j(t^2)) \o_R u_j \\
\cop_T(t) & := &\sum_i t_i \o_R (\beta_i(t^1) \o_B t^2)
\end{eqnarray}
and counit $\eps_T(t) = t^1 t^2$, which is the multiplication mapping $\mu$ restricted to $T$.
It is also true for $T$ that it is a right $R$-bialgebroid when $A \| B$ is in possession of only a left D2 quasibase,
or only a right D2 quasibase; the only possible loss being finite $R$-projectivity of $T$ on one side.  

We note that $S$ and $T$ are $R$-dual bialgebroids w.r.t.\ either of the nondegenerate pairings $\bra \alpha \| t \ket := \alpha(t^1)t^2$ or $[ \alpha \| t ] := t^1 \alpha(t^2) $, both
with values in $R$ \cite[5.3]{KS}. 
 
The opposite multiplication $\overline{\mu}$ on $A\op$ is characterized by the equation
$\overline{\mu}(\overline{a} \o \overline{b}) = \overline{\mu(b \o a)}$ using the notation in Prop.~\ref{prop-op}.
  The next proposition shows that passing to opposite algebras does not change the chirality of the bialgebroids
$S$ and $T$
 associated with a D2 extension, but it does change their underlying $R$-corings to the co-opposites.  

\begin{prop}
\label{prop-teecop}
Let $A \| B$ be an algebra extension.  Then $A \| B$ is left D2 with $R$-bialgebroids $S$
and $T$ if and only if $A\op \| B\op$ is right D2 with $R\op$-bialgebroids $S_{\rm cop}$
and $T_{\rm cop}$.
\end{prop}
\begin{proof}
Suppose $A \| B$ is left D2 with $a \o a' = \sum_i t_i \beta_i(a) a'$.  Then
we have seen in Prop.~\ref{prop-op} that $A\op \| B\op$ is right D2; in fact,
$\overline{a \o a'} = \sum_i \overline{a'}\, \overline{\beta_i}(\overline{a}) \overline{t_i}$
is the equation of a right D2 quasibase. The centralizer $C_{A\op}(B\op) = R\op$ is the
opposite of $R = C_A(B)$.  
If we denote the comultiplication in $T(A\op \| B\op)$ by $\overline{\cop}_T$ we have 
$$ \overline{\cop}_T(\overline{t}) = \sum_i \overline{t^2} \o \overline{\beta_i}(\overline{t^1}) \o_{R\op} \overline{t^2_i} \o \overline{t^1_i} = \overline{\cop_T(t)} $$
We conclude that $\overline{\cop}_T = \cop^{\rm op}_T$.
Similarly $\overline{\cop}_S = \cop^{\rm op}_S$.  It is also apparent that source and target maps are reversed
in passing to the opposite algebras, e.g., $\overline{\lambda(a)} = = \rho(\overline{a})$
for $a \in A$.  Finally we have shown in Prop.~\ref{prop-op} that $S( A\op \| B\op) \cong S$
and $T(A\op \| B\op) \cong T$ as algebras. The converse is proven similarly.   
\end{proof}

In \cite[4.1]{KS} we observed that $S$ acts on $A$ 
via evaluation as a left $S$-module algebra (or algebroid): if $A_B$ is a balanced module,
then the invariant subalgebra $A^S = B$.  In this paper, we will be more concerned with the dual
concept, comodule algebra (defined below).  As an example of this duality, and a guide
to what we are about to do, we dualize, as we would (but more carefully)
 for  Hopf algebra actions, 
the left action just mentioned $\lact: S \o_R A \to A$,  $\alpha \lact a := \alpha(a)$ for $\alpha \in S, a \in A$, to a right coaction
$\varrho_T: A \to A \o_R T$ given by $\rho(a) = a\0 \o a\1$ where
$\alpha \lact a = $ $ a\0 [\, \alpha \, |\, a\1 \, ]$. This comes out as $\varrho_T(a) = 
\sum_j \gamma_j(a) \o u_j$, since $\alpha(a) =$ $ \sum_j \gamma_j(a) [\alpha \| u_j]$  (obtained by applying
$\id \o \alpha$ to eq.~(\ref{eq: right d2 qb})). The resulting right $T$-comodule algebra structure on $A$
is studied in \cite{fer}, where it is shown that assuming $A_B$ balanced and $A \| B$ right D2
results in a Galois extension $A \| B$ 
in the usual Galois coaction picture.

There is also an action of $T$ on $\ME$ studied 
in \cite[5.2]{KS}: the $R$-bialgebroid $T$ acts from the right on $\ME$ by
\begin{equation}
\label{eq: ract by tee}
f \ract t := t^1 f(t^2 -)
\end{equation} for $f \in \ME, t \in T$.
This action makes $\ME$ a right $T$-module algebra
with invariants $\rho(A)$ (where recall $\rho(a)(x) = xa$ for $x,a \in A$). Thinking in terms of Hopf algebra
duality, we then expect to see a left coaction $\varrho: \ME \to S \o_R \ME$ with Sweedler
notation $\varrho(f) = f\-1 \o f\0$ 
satisfying $f \ract t = [\, f\-1 \, | \, t \, ] f\0$. This comes out in terms of a right D2 quasibase as 
\begin{equation}
\label{eq: row}
\varrho(f) = \sum_j \gamma_j \o (f \ract u_j)
\end{equation}
since $f \ract t = $ $\sum_j t^1 \gamma_j(t^2)u_j^1 f(u_j^2-) = $
$ \sum_j [\, \gamma_j \, |\, t \, ] (f \ract u_j)$.  Since
$\ME$ is a variant of a smash product of $A$ with $S$
(cf.\ \cite[section 3]{KS}), we would want to 
show that $\ME$ is a Galois extension of a copy of $A^{\rm op}$
somewhat in analogy with cleft Hopf algebra coaction
theory although there is no antipode in our set-up:  see eq.~(\ref{eq: total integral}) for why
we may think of the natural inclusion $S \into \ME$ as
a total integral which cleaves the $S$-extension $\ME \| \rho(A)$.    
We next turn to several definitions,
lemmas and the main theorem below in which we prove that $\varrho$ is a Galois coaction
for the left $S$-extension $\ME$ over $\rho(A)$.

\begin{definition}
Let $S'$ be a left $R'$-bialgebroid $(S', \tilde{s}, \tilde{t},$ $ \cop, $ $ \eps)$.   
A left $S'$-comodule algebra $C$ is an algebra $C$ with algebra homomorphism $R' \to C$ 
 together with a coaction $\delta: C \to S' \o_{R'} C$, where  values $\delta(c)$ are denoted by the Sweedler
notation $c\-1 \o c\0$, 
such that $C$ is a left $S'$-comodule over the $R'$-coring $S'$ \cite[18.1]{BW}, 
\begin{equation}
 \delta(1_C) = 1_{S'} \o 1_C, 
\end{equation} 
\begin{equation}
\label{eq: homog}
c\-1 \tilde{t}(r) \o c\0 = c\-1 \o c\0 \cdot r
\end{equation}
for all $r \in R$,
and 
\begin{equation}
\label{eq: cop homo}
\delta(cc') = \delta(c) \delta(c')
\end{equation}
 for all $c,c' \in C$.   The subalgebra of coinvariants
is \begin{equation}
{}^{\rm co \, S'}C := \{ c \in C \| \delta(c) = 1_{S'} \o c \}.
\end{equation} 
$C$ is said to be a left $S'$-extension of ${}^{\rm co \, S'}C $. 
\end{definition} 

Like in the definition of left bialgebroid, the axiom~(\ref{eq: cop homo}) makes sense
because of the axiom~(\ref{eq: homog}). 
The algebra homomorphism $R' \to C$ induces a natural $R'$-bimodule
on $C$ which we refer to implicitly.

\begin{example}
\begin{rm}
The left bialgebroid $S'$ is a left $S'$-comodule algebra with coinvariants $\Im \tilde{t}$. 
For let the coaction coincide with the comultiplication, the algebra homomorphism of the base algebra $R'$
into the total algebra $S'$ be the source mapping $\tilde{s}$, so that eq.~(\ref{eq: homog}) follows
from eq.~(\ref{eq: duo}) and computing with $r \in R'$:
$$ \cop (\tilde{t}(r)) = \cop(1 \cdot r) = \cop(1) \cdot r = 1 \o_{R'} \tilde{t}(r) $$
so $\tilde{t}(R') \subseteq {}^{\rm co \, S'}S' $.  Conversely, if $\cop(c) = 1 \o c$, then
$$ c= c\1 \cdot \eps(c\2) = 1\cdot \eps(c) = \tilde{t}(\eps(c)) $$
so that $c \in \Im \tilde{t}$, whence $\tilde{t}(R') = {}^{\rm co \, S'}S' $.
\end{rm}
\end{example}

\begin{definition}
Let $S'$ be a  left $R'$-bialgebroid.  
An $S'$-comodule algebra $C$ is a left $S'$-Galois extension of its coinvariants
$D$ if the (Galois) mapping 
$\beta: C \o_D C \to S' \o_{R'} C$ defined
by $\beta(c \o c') = c\-1 \o c\0 c' $ is bijective.
\end{definition} 

The definition of Galois extension is equivalent to $S' \o_{R'} C$ being a Galois $C$-coring \cite{BW}.
The definitions of right comodule algebra
and right Galois extension for right bialgebroids are very similar
(and found in \cite[1.2, 1.4]{fer}, cf.\  \cite[31.23]{BW}).  
The two notions are in fact in one-to-one correspondence via opposite algebras as we show next.

\begin{prop}
\label{prop-cop}
Let $A \| B$ be an algebra extension and $T$ a right bialgebroid over an algebra $R$.
Then $A \| B$ is a right $T$-extension (or Galois extension) if and only if
$A\op \| B\op$ is left $T^{\rm op}_{\rm cop}$-extension (or Galois extension, respectively).
\end{prop}
\begin{proof}
It is noted in \cite{KS}, and explored further in the next subsection in this paper,
that $T =$ $(T, R, s_R, t_R, \cop, \eps)$ is a right $R$-bialgebroid iff $T^{\rm op}_{\rm cop}=$$(T\op, R\op, s_R, t_R,
\cop\op, \eps)$ is a left $R\op$-bialgebroid.  
If $A \| B$ is a right $T$-extension with coaction $\delta(a) = a\0 \o a\1$, then we may define a $T^{\rm op}_{\rm cop}$-algebra structure on $A\op$
via the induced homomorphism $R\op \to A\op$,
$$\overline{\delta}(\overline{a}) = \overline{a}\-1 \o \overline{a}\0 := \overline{a\0 \o a\1}$$ using an obvious extension of the notation in Prop.~\ref{prop-op}. All the properties of a left comodule algebra
 follow routinely, including $A^{\rm co \, T} = {}^{\rm co \, T\op_{\rm cop}}A$.  
The one-to-one correspondence of right and left Galois extensions follows
from noting that the Galois mapping $\beta: A\op \o_{B\op} A\op \to
T^{\rm op}_{\rm cop} \o_{R\op} A\op$ defined above and the right Galois mapping $\beta_R: A \o_B A \to A \o_R T$
given by 
$\beta_R(a \o a') := a {a'}\0 \o {a'}\1$ are isomorphically related by a commutative square, i.e., 
$$\beta(\overline{a \o_B a'}) = \overline{\beta_R(a \o_B a')} $$
for $a, a' \in A$.  
\end{proof}

\subsection*{A new example of Hopf algebroid}
One unfortunate side-effect of generalizing Lu's bialgebroid $A^e$ over an
algebra $A$ \cite{Lu} to the right bialgebroid $T$ of a D2 extension
$A \| B$ (in \cite{KS}) is that the antipode is lost, for the flip or twist anti-automorphism on $A \o_K A^{\rm op}$ 
does not extend to a self-mapping of $(A \o_B A)^B$. If we require $B$ to be separable with symmetric separability
element however, there is a projection of $A \o_B A \to T$ which we may apply to define
a twist of $T$.  However, 
Lu's definition of Hopf algebroid \cite{Lu} makes it necessary to find an appropriate section $T \o_R T \to T \o_K T$ of the canonical
map in the other direction, although  the centralizer $R$
is not \textit{a priori} separable.  In this subsection, we carry out 
this plan using instead the alternative definition
of Hopf algebroid in \cite[B\"{o}hm and Szlach\'{a}nyi]{BS}. 
They define a Hopf algebroid to be a left $R'$-bialgebroid $S'$ with anti-automorphism $\tau: S' \to S'$ (called antipode)
   satisfying eqs.~(\ref{eq: (4.1)}),~(\ref{eq: (4.2)}) and~(\ref{eq: (4.3)})
below.  
 
Let $K$ be a commutative ring and $B$ a \textit{Kanzaki separable $K$-algebra}
\cite[strongly separable algebra]{KSt}.  This means that there is a separability element
$e = e^1 \o e^2$ $ \in B \o_K B$ which is symmetric, so
that $e^1 e^2 = 1 $ $ = e^2 e^1$ as well as $be = eb$
and $e^1 b \o e^2 = $ $ e^1 \o be^2$ for all $b \in B$. (Typically
for quantum algebra, we use both
these equalities repeatedly below together with $bt = tb$ for $t \in T$ and well-definedness of various mappings on equal commuting elements.)
For example, all separable algebras over a field of characteristic zero
are Kanzaki separable.  Over a field of characteristic $p$ matrix
algebras of order divisible by $p$ are separable although not
Kanzaki separable.  Fix the notation above for the next theorem.

\begin{theorem}
\label{th-newex}
Let $A \| B$ be a D2 extension of $K$-algebras where
$B$ is Kanzaki separable. Then the left bialgebroid
$T^{\rm op}_{\rm cop}$ is a Hopf algebroid.
\end{theorem}
\begin{proof}
The standard right bialgebroid $T = (A \o_B A)^B$ with structure $(T, R, s_R, t_R,$ $ \cop, \eps)$ becomes a left
bialgebroid via the opposite multiplication and co-opposite comultiplication as follows \cite[2.1]{KS}:  
$$T^{\rm op}_{\rm cop} := ( T^{\rm op}, R^{\rm op}, s_L = s_R, t_L = t_R, \cop^{\rm op}, \eps). $$   
The multiplication on $T^{\rm op}$ is 
\begin{equation}
tt' = t^1 {t'}^1 \o_B {t'}^2 t^2 
\end{equation}
while the target and source maps are $s_L: R^{\rm op} \to T^{\rm op}$, $t_L: R \to T^{\rm op}$ are then
$s_L(r) = 1_A \o r$, $t_L(r) = r \o 1_A$ for $r \in R\op$.  The $R^{\rm op} $ - $R^{\rm op}$-bimodule structure
on $T$ is then given by
\begin{equation}
r \cdot t \cdot r' = s_L(r) t_L(r') t = (r' \o r)t = r' t^1 \o t^2 r. 
\end{equation}
In other words ${}_{R^{\rm op}}T_{R^{\rm op}}$ is the standard bimodule ${}_RT_R$ with endpoint multiplication after passing
to modules over the opposite algebra of $R$.  Tensors over $R^{\rm op}$ are the same as tensors over $R$ after a flip;
e.g., $$
T \o_{R^{\rm op}}T \stackrel{\cong}{\longrightarrow} T \o_R T \stackrel{\cong}{\longrightarrow} (A \o_B A \o_B A)^B 
$$
via the mapping 
\begin{equation}
\label{eq: iq}
t \o t' \stackrel{\cong}{\longmapsto} {t'}^1 \o {t'}^2 t^1 \o t^2,
\end{equation}
which is an $R^{\rm op}$-$R^{\rm op}$-isomorphism \cite[5.1]{KS}. 

The comultiplication $\cop^{\rm op}: T \to T \o_{R^{\rm op}} T$ is given for $t \in T$ by
the co-opposite of eq.~(\ref{eq: tee}), 
\begin{equation}
\label{eq: tee cop}
 \cop^{\rm op}(t)  = \sum_j u_j \o (t^1 \o_B \gamma_j(t^2))
\end{equation}
which in $(A \o_B A \o_B A)^B$ is the value $t^1 \o 1_A \o t^2$ after applying eq.~(\ref{eq: iq})
and the right D2 quasibases eq.~(\ref{eq: right d2 qb}). We will denote below $\cop^{\rm op}(t) = $ $t\1 \o t\2$. 

The antipode $\tau: T \to T$ is a flip  composed with a projection from $A \o_K A$:
\begin{equation}
\label{eq: antipode}
\tau(t) := e^1 t^2 \o_B t^1 e^2.
\end{equation}
 We next note that $\tau$ is an anti-automorphism of order two on $T^{\rm op}$. Let $e = f$ in $B \o_K B$ so that  
$$ \tau(t') \tau(t) = e^1 {t'}^2 f^1 t^2 \o_B t^1 f^2 {t'}^1 e^2 = e^1 {t'}^2 t^2 \o_B t^1 f^2 f^1 {t'}^1 e^2 = $$
$ = \tau(tt')$ since ${t'}^1 e^2 \o_K e^1 {t'}^2 \in (A \o_K A)^B$ and $f^2 f^1 = 1_B$. In addition,
$$ \tau^2(t) = \tau(e^1 t^2 \o t^1 e^2) = f^1 t^1 e^2 \o_B e^1 t^2 f^2 = t, $$
since $e^2 e^1 = 1$.

Next, we show that $\tau$ satisfies the three axioms of \cite[Def.\ 4.1]{BS} given below in
eqs.~(\ref{eq: (4.1)})-(\ref{eq: (4.3)}).  Note that $\tau^{-1} = \tau$.  First,
\begin{equation}
\label{eq: (4.1)}
\tau \circ t_L = s_L
\end{equation}
since $\tau(r \o 1) = e^1 \o_B r e^2 = 1 \o r$ for $r$ in the
centralizer $C_A(B)$.  

Second, we have the equality in $T\o_{R^{\rm op}} T$, 
\begin{equation}
\label{eq: (4.2)}
\tau^{-1}(t\2)\1 \o \tau^{-1}(t\2)\2 t\1 = \tau^{-1}(t) \o 1_T
\end{equation}
which follows from eq.~(\ref{eq: tee cop}) and (working from left to right): 
$$ \sum_{j,k} u_k \o (e^1 \gamma_j(t^2)  \o_B  \gamma_k(t^1 e^2))u_j =
\sum_{j,k} u_k \o (e^1 \gamma_j(t^2)u_j^1 \o_B  u_j^2\gamma_k(t^1 e^2)) $$
the last expression mapping via the isomorphism in eq.~(\ref{eq: iq}) into
$$e^1 \o_B t^2 \gamma_k(t^1) u_k^1 \o_B u^2_k e^2 = 1_A \o e^1 t^2 \o t^1 e^2 $$
in $(A \o_B A \o_B A)^B$, which is the image of $\tau^{-1}(t) \o 1_T$ under the
same isomorphism.  

Finally, we have the equality in $T\o_{R^{\rm op}} T$,
\begin{equation}
\label{eq: (4.3)}
\tau(t\1)\1 t\2 \o \tau(t\1)\2 = 1_T \o \tau(t)
\end{equation}
for all $t \in T$, which follows similarly from
$$\sum_{j,k} u_k(t^1 \o_B \gamma_j(t^2)) \o e^1 u^2_j \o_B \gamma_k(u^1_j)e^2) = 
\sum_{j,k} u_k^1 t^1 \o_B \gamma_j(t^2)u_k^2) \o e^1 u^2_j \o_B \gamma_k(u^1_j)e^2)$$
which maps isomorphically to
$$\sum_{j,k} e^1 u^2_j \o_B \gamma_k(u^1_j)u^1_k t^1 \o_B \gamma_j(t^2) u^2_k e^2 = \sum_j
e^1 u^2_j \o t^1 \o \gamma_j(t^2) u^1_j e^2 = $$
$$ e^1 t^2 \o t^1 e^2 \o 1_A \stackrel{\cong}{\longleftarrow} 1_T \o \tau(t) $$
since $\sum_j \gamma_j(a)u^1_j e^2 \o_K e^1 u^2_j = e^2 \o_K e^1a$ for $a \in A$
follows from eq.~(\ref{eq: right d2 qb}). 
\end{proof}
  
Suppose $K$ is a field, then $B$ is finite dimensional as it is a separable algebra. If moreover 
 $A \| B$ is a proper algebra extension, then there is bimodule projection $A \to B$ by separability,
whence $A_B$ is finitely generated (and projective by the D2 quasibase eq.~(\ref{eq: right d2 qb})) and so $A$ is finite dimensional as well.  
The theorem is thus viewed as a natural generalization of Lu's Hopf algebroid with twist
antipode to certain finite-dimensional algebra pairs.


\section{Main Theorem: endomorphism algebra extension is Galois}
 
Consider a right D2 extension $A \| B$ with right D2 quasibases $\gamma_j \in S$,
$u_j \in T$, left endomorphism algebra $\ME = \End {}_BA$
and natural algebra embedding $\rho: A \into \ME$ by right multiplications.  We
introduced in eq.~(\ref{eq: row}) a possible coaction $\varrho: \ME \to S \o_R \ME$, which in this section we show  
 makes the extension $\ME \| \rho(A)$ into a left Galois extension --- a 
 nontrivial proof involving several lemmas at first.
The next lemma will be used among other things to show that the coaction
$\varrho$ is coassociative. 

\begin{lemma}
\label{lemma: a}
Let $A \| B$ be a right D2 extension. Then we have the isomorphisms
\begin{equation}
\label{eq: short iso}
S \o_R \ME \cong \Hom ({}_BA \o_B A, {}_B A) 
\end{equation}
via $\alpha \o f \longmapsto (a \o a' \mapsto \alpha(a)f(a'))$, and
\begin{equation}
\label{eq: long iso}
S \o_R S \o_R \ME \cong \Hom ({}_BA \o_B A \o_B A, {}_B A)
\end{equation}
via $\alpha \o \beta \o f \longmapsto (a \o a' \o a'' \mapsto
\alpha(a) \beta(a') f(a''))$.
\end{lemma}
\begin{proof}
The inverse in (\ref{eq: short iso}) is given by $F \mapsto $
$ \sum_j \gamma_j \o u^1_j F(u^2_j \o -) $ by eq.~\ref{eq: right d2 qb}.

The inverse in (\ref{eq: long iso}) is given by $$F \mapsto 
\sum_{j,k,i} \gamma_j \o \gamma_k \o u^1_k u^1_j F(u^2_j \gamma_i(u^2_k)u^1_i \o u^2_i \o -) 
$$ since
$$ \sum_{j,k,i} \gamma_j(a)\gamma_k(a')u^1_k u^1_j F(u^2_j \gamma_i(u^2_k)u^1_i \o u^2_i \o a'') 
 = \sum_{j,k} \gamma_j(a)u^1_jF(u^2_j \gamma_i(a') u^1_i \o u^2_i \o a'')$$
$ = F(a \o a' \o a'') $
and
$$ \sum_{j,k,i} \gamma_j \o \gamma_k \o_R u^1_k u_j^1 \alpha(u^2_j \gamma_i(u^2_k)u^1_i) \beta(u^2_i) f(-) =
\sum_{j,i} \gamma_j \o u^1_j \alpha(u^2_j \gamma_i(-)u^1_i) \beta(u^2_i) \o f(-) $$
$ = \sum_j \gamma_j(-) u^1_j \alpha(u^2_j) \o \beta \o f = \alpha \o \beta \o f$
for $f \in \ME, \alpha, \beta \in S$. 
\end{proof}

The existence alone of an isomorphism in the next lemma may
be seen by letting $M$ be free of rank one. 

\begin{lemma}
Given any algebras $A$ and $B$, with modules $M_A$, ${}_BU$ and bimodule ${}_BN_A$ with $M_A$ f.g.\ projective,
then 
\begin{equation}
\label{eq: mull}
M \o_A \Hom({}_B N, {}_B U) \stackrel{\cong}{\longrightarrow} \Hom ({}_B\Hom (M_A, N_A), {}_B U) 
\end{equation}
via the mapping $m \o \phi \longmapsto (\nu \mapsto \phi (\nu(m)))$.
\end{lemma}
\begin{proof}
Let $m_i \in M$, $g_i \in \Hom (M_A, A_A)$ be dual bases for $M_A$.
For each $n \in N$ let $ng_i(-)$ denote the obvious mapping
in $\Hom (M_A, N_A)$.  
Then the inverse mapping is given by 
\begin{equation}
F \longmapsto \sum_i m_i \o (n \mapsto F(ng_i(-))).
\end{equation}
  Note that both maps are well-defined module homomorphisms,
and inverse to one another since 
$$ \sum_i m_i \o \phi(-g_i(m)) = \sum_i m_i \o_A (g_i(m) \phi)(-) =
\sum_i m_i g_i(m) \o \phi = m \o \phi$$
and $ \nu \mapsto \sum_i F(\nu(m_i) g_i ) = F(\nu)$ for
$\nu \in $ $\Hom (M_A, N_A)$. 
\end{proof}

The lemma above is relevant in our situation since the D2 condition
implies that a number of constructions such as the tensor-square and 
endomorphism algebras are finite projective. For example, $\ME_A$ is f.g.\ 
projective \cite[3.13]{KS}, which we may also see directly from eq.~(\ref{eq: right d2 qb})
by applying $\id_A \o_B f$ for $f \in \ME$, viewing $\gamma_j \in S \subseteq \ME$
and an obvious mapping  of $A \o_B A$ into $\Hom (\ME_A, A_A)$ which appears in the next
lemma.

\begin{lemma}
\label{lemma: b}
If $A \| B$ is right D2, then
$${}_B A \o_B A \stackrel{\cong}{\longrightarrow} {}_B \Hom(\ME_A, A_A) $$
via $\Psi(a \o a')(f) = af(a') $.
\end{lemma}
\begin{proof}
 Let $F \in (\ME_A)^*$ (the right $A$-dual of $\ME$).  Define an inverse 
$\Psi^{-1}(F) = \sum_j  F(\gamma_j)u_j$.  Then $\Psi^{-1} \Psi = \id_{A \o_B A}$ by eq.~(\ref{eq: right d2 qb}).
Also $\Psi \Psi^{-1} = \id_{\ME^*}$ since for $f \in \ME$, $$\Psi \Psi^{-1}(F)(f) = \sum_j F(\gamma_j)u^1_j f(u^2_j) = F(\sum_j \gamma_j(-) u^1_j f(u^2_j)) = F(f). \qed $$ 
\renewcommand{\qed}{}\end{proof}

Recall that $\rho(A)$ denotes the set of all right multiplication operators by elements of $A$.  

\begin{theorem}
\label{th-endogalext}
Let $A \| B$ be a right D2 extension.  Then $\ME$ is a left $S$-comodule algebra with the coaction~(\ref{eq: row})
and a Galois extension of its coinvariants $\rho(A)$.  
\end{theorem}

\begin{proof}
The coaction $\varrho$ in eq.~(\ref{eq: row}) has value on $f \in \ME = \End {}_BA$  given by $$ f\-1 \o f\0 = \sum_j \gamma_j \o u_j^1 f(u_j^2 -)$$
where $\gamma_j \in S, u_j \in T$ is a right D2 quasibase. First, the algebra homomorphism $R \to \ME$ is 
given by $\lambda$, so 
for $1_{\ME} = \id_A = 1_S$, we have 
$$ \varrho(\id_A) = \sum_j \gamma_j \o_R \lambda(u^1_j u^2_j) = \sum_j \gamma_j(-)u_j^1 u_j^2 \o \id_A = 1_S \o 1_{\ME}.
$$

Secondly, we show that $\ME$ forms a left $S$-comodule w.r.t.\ the $R$-coring $S$ and the coaction $\varrho$. 
The coaction is coassociative, $(\cop_S \o \id_{\ME})\varrho =$ $ (\id_S \o \varrho)\varrho$,
for we use lemma~\ref{lemma: a} (as an identification and suppressing the isomorphism) and eqs.~(\ref{eq: right d2 comult})
and~(\ref{eq: tee mult})  to check values of each side of this equation, evaluated on $A \o_B A \o_B A$:
\begin{eqnarray*}
 (\sum_j \cop_S(\gamma_j) \o (f \ract u_j))(a \o a' \o a'') &=& \sum_{k,j} \gamma_k(a)u_k^1\gamma_j(u^2_k a')u_j^1 f(u_j^2 a'')  \\
& = & \sum_j \gamma_j(aa') u_j^1 f(u^2_ja'') = f(aa'a'') \\
& = & \sum_{i,j} \gamma_i(a) \gamma_j(a')u^1_j u^1_i f(u^2_i u^2_j a'') \\
& = & (\sum_j \gamma_j \o \varrho (f \ract u_j))(a \o a' \o a''). 
\end{eqnarray*}
We  also show that $\varrho: \ME \to S \o_R \ME$ is a left $R$-module morphism: given $r \in R, f\in \ME$, we
 use lemma~\ref{lemma: a} again to note that for $a, a' \in A$ 
$$ \varrho(\lambda(r)f) (a \o a') = \sum_j \gamma_j(a)u^1_j rf(u^2_j a') = rf(aa') = (r\cdot f\-1 \o f\0)(a \o a'), $$
by an application of eq.~(\ref{eq: right d2 qb}) (inserting an $r \in C_A(B)$).  
(Note with $r = 1$ that we obtain  
\begin{equation}
\label{eq: total integral}
\varrho(\alpha) = \cop_S(\alpha) \ \ \ \ \ \forall \, \alpha \in S,
\end{equation}
which should be compared to the total integral and cleft extension approach in \cite[6.1]{KN} and \cite{DT}.)
Finally, $\ME$ is counital, whence a left $S$-comodule, since for $f \in \ME$, 
$$ (\eps_S \o \id_S)\varrho(f) = \sum_j \gamma_j(1)u^1_jf(u^2_j -) = f. $$

Next we compute that $\Im \varrho$ lies in a submodule of $S \o_R \ME$ where tensor product multiplication
makes sense: again using lemma~\ref{lemma: a} and for $a,a' \in A$,
$$ (f\-1 \tilde{t}(r) \o f\0)(a \o a') = \sum_j \gamma_j(ar)u_j^1 f(u_j^2 a') = f(ara') $$
$$ = \sum_j \gamma_j(a) u_j^1 f(u_j^2 ra') = (f\-1 \o f\0 \lambda(r))(a \o a'). $$
Then multiplicativity of the coaction follows from the measuring axiom satisfied by the right action of
$T$ on $\ME$ \cite[5.2]{KS} and eq.~(\ref{eq: tee}) ($f,g \in \ME$):  
$$ \varrho(fg)(a \o a') = \sum_j (\gamma_j \o (f\ract {u_j}\1) \circ (g \ract {u_j}\2))(a \o a') = $$
$$ \sum_{j,k} \gamma_j(a) u^1_j f(\gamma_k(u^2_j) u^1_k g(u^2_k a')) = f(g(aa')) = \sum_{j,i} \gamma_i(\gamma_j(a))
u^1_i f(u^2_i u^1_j g(u^2_j a')) = $$
$$ = \sum_{i,j} (\gamma_i \circ \gamma_j \o (f \ract u_i ) \circ (g \ract u_j))(a \o a') = \varrho(f)\varrho(g)(a \o a'). $$

Next we determine the coinvariants ${}^{\rm co \, S} \ME$. Given $a \in A$,
we note that $\rho(a) \in {}^{\rm co \, S} \ME$ since
$$ \varrho(\rho(a)) = \sum_j \gamma_j \o u^1_j u^2_j(-a) = 1_S \o \rho(a).$$
Conversely, suppose $\sum_j \gamma_j \o (f \ract u_j) = 1_S \o f$
in $S \o_R \ME \cong \Hom ({}_BA \o_B A, {}_BA)$,
then for $a, a' \in A$,
$$ \sum_j \gamma_j(a)u^1_j f(u^2_j a') = f(aa') = af(a'). $$
It follows that $f(a) = af(1_A)$ for all $a \in A$, so
$f = \rho(f(1)) \in \rho(A)$.  Hence, ${}^{\rm co \, S} \ME = \rho(A)$. 

Finally, the Galois mapping 
\begin{equation}
\label{eq: gm}
\beta: \ME \o_{\rho(A)} \ME \to
S \o_R \ME, \ \ \ \ \ \beta(f \o g) = f\-1 \o f\0 g 
\end{equation}
under the identification $S \o_R \ME \cong$ 
$\Hom ({}_B A \o_B A, {}_BA)$ in lemma~\ref{lemma: a}
is given by an application of eq.~(\ref{eq: right d2 qb}):  ($a,a' \in A, f,g \in \ME$)
 
\begin{equation}
\label{eq: galois map}
\beta(f \o g)(a \o a') = \sum_j \gamma_j (a)u^1_j f(u^2_j g(a')) = f(ag(a')).
\end{equation}

We show $\beta$ to be a  composite of several isomorphisms using the lemmas
(and see the commutative diagram in section~1).  First note
that 
$\rho(A) \cong  A\op$ 
and 
${}_{\rho(A)} 
\ME_{\rho(A)}$
 given by
 $\rho(a')\circ f \circ \rho(a) (a'') = f(a'' a) a'$ 
is equivalent to ${}_A \ME_A$ given by $a \cdot f \cdot a' (a'') = f(a'' a)a'$.  
This is the usual $A$-$A$-bimodule structure on the left endomorphism algebra $\ME$
considered in \cite[3.13]{KS}, where  $\ME_A$ is shown to be f.g.\ projective.
Consider then the composition of isomorphisms,
\begin{eqnarray*}
 \ME \o_{\rho(A)} \ME  & \stackrel{\cong}{\longrightarrow} & \ME \o_A \Hom ({}_BA, {}_BA) \\
& \stackrel{\cong}{\longrightarrow} & \Hom ({}_B \Hom(\ME_A, A_A), {}_BA) \\
& \stackrel{\cong}{\longrightarrow} & \Hom ({}_B A \o_B A, {}_B A) 
\end{eqnarray*}
given by $$ f \o g \mapsto g \o f \mapsto (\nu \mapsto f(\nu(g))) \mapsto ( a\o a' \mapsto  f( a g(a')). $$
This is $\beta$ as given in eq.~(\ref{eq: galois map}), whence 
$\beta$ is an isomorphism and the extension $\ME \| \rho(A)$ is  Galois.  
\end{proof}

By chasing the diagram in Figure~1 around in the opposite direction,
we obtain the inverse Galois mapping, $\beta^{-1}: S \o_R \ME \to \ME \o_{\rho(A)} \ME$,
given by 
\begin{equation}
\label{eq: beta inverse}
\beta^{-1}(\alpha \o h) = \sum_j \alpha(- u^1_j) h(u^2_j) \o \gamma_j.
\end{equation}
We check this directly ($f,g,h \in \ME$, $\alpha \in S$):
$$\beta^{-1}(\beta(f \o g)) = \beta^{-1}(\sum_j \gamma_j \o u^1_j f(u^2_j g(-)) = 
\sum_{j,k} \gamma_j(-u^1_k)u^1_j f(u^2_j g(u^2_k)) \o \delta_k
= $$
$$ \sum_k f(-u^1_k g(u^2_k)) \o_{\rho(A)} \gamma_k = \sum_k f(-) \o \rho(u^1_kg(u^2_k)) \circ \gamma_k = f \o g,
$$
by eq.~(\ref{eq: right d2 qb}) and $$ \beta(\beta^{-1}(\alpha \o h) = \sum_{j,k} \gamma_j \o u^1_j \alpha(u^2_j\gamma_k(-)u^1_k)h(u^2_k) =
\sum_j \gamma_j \o_R u^1_j \alpha(u^2_j) h(-) = \alpha \o h. $$

\section{Corollaries and a converse of the main theorem}

In this section we recall the characterization in \cite{fer} of Galois extension
as one-sided D2 and balanced extensions.  Then 
 right D2 extensions are shown to have left 
endomorphism algebra extensions that are left D2.
Right
endomorphism algebra extensions are shown to be right D2. A converse endomorphism ring theorem for one-sided split D2 Frobenius
extensions is proven in Theorem~{th-converse-endo}.

We next recall the right D2 characterization of a right Galois $T$-extension
$A \| B$ for a right bialgebroid $T$ over an algebra $R$
\cite{fer} (with roots in \cite{KN, NV}). We show in detail why this is equivalent
to a left D2 characterization of a left Galois extension.

\begin{theorem}[\cite{fer}]
\label{th-characterization}
Let $A \| B$ be a proper algebra extension. Then the following hold
and are equivalent: 
\begin{enumerate}
\item \label{item-en} 
$A \| B$ is a right $T$-Galois extension for
some left finite projective right bialgebroid $T$ over
some algebra $R$ if and only if $A \| B$ is right D2 and right balanced.
\item \label{item-to}
$A \| B$ is a left $S$-Galois extension for
some right finite projective left bialgebroid $S$ over
some algebra $R$ if and only if $A \| B$ is left D2 and left balanced.
\end{enumerate}
\end{theorem}
\begin{proof}
We have established the characterizations in (\ref{item-en}) or
(\ref{item-to}) in \cite{fer}. Their equivalence may be seen
as follows.  Suppose we assume (\ref{item-en}) and we are
given a left D2 and left balanced extenson $A \| B$.
Then by Prop.~\ref{prop-op} $A\op \| B\op$ is right D2,
and moreover right balanced since $\End {A\op}_{B\op}$
$= \End {}_BA$ in which  left multiplication $\lambda(\overline{a})$
is equal to $\overline{\rho(a)}$, the canonical image of right multiplication. Then $A\op \| B\op$ is a right $T_{\rm cop}$-Galois
extension where $T_{\rm cop}$ is the associated right $R^{\rm op}$-bialgebroid by Prop.~\ref{prop-teecop}, whence by Prop.~\ref{prop-cop}
 $A \| B$ is a left $T\op$-Galois extension
for the left $R$-bialgebroid $T^{\rm op}$. Since $T_R$ is f.g.\
projective by the left D2 hypothesis, it follows that
$T^{\rm op}$ is right f.g.\ $R$-projective as well,
since $T$ and $T\op$ have the same underlying $R$-$R$-bimodule.

The converse follows by reversing the order of application 
of Props.~\ref{prop-op},~\ref{prop-teecop} and~\ref{prop-cop}, and noting that $S^{\rm op}_{\rm cop}$ is left finite projective over
$R\op$ 
iff $S$ is right finite projective over $R$.  
\end{proof}

The following endomorphism ring theorem is then a corollary of this theorem and the main Theorem~\ref{th-endogalext}.

\begin{cor}[Endomorphism Ring Theorem for D2 Extensions]
\label{cor-leftendo}
If $A \| B$ is a right D2 algebra extension, then $\ME \| A\op$
is a left D2 extension (and left balanced).  
\end{cor}

This generalizes \cite[Theorem 6.7]{KS} by
leaving out the hypothesis that $A \| B$ be a Frobenius extension
(i.e., $A_B$ be finite projective and $A \cong \Hom (A_B, B_B)$
as natural $B$-$A$-bimodules). Via eq.~(\ref{eq: beta inverse}) one 
arrives at left D2 quasibases for $\ME \| A\op$ in terms of right D2 quasibases for
$A \| B$:  

\begin{equation}
\label{eq: left D2 qb endo}
 T_j := \sum_k \gamma_j (- u_k^1) u^2_k \o \gamma_k, \ \ \ \ \mathcal{B}_j(f) := f \ract u_j = u_j^1 f(u_j^2 -) 
\end{equation}
We note that $T_j = \beta^{-1}(\gamma_j \o \id_A) \in (\ME \o_{\rho(A)} \ME)^{\rho(A)}$ and $\mathcal{B}_j \in \End {}_{\rho(A)}
\ME_{\rho(A)}$, and verify the essential part of the left D2
quasibase eq.~(\ref{eq: left d2 qb}) in $\ME \o_{\rho(A)} \ME$: ($f \in \ME$)

$$\sum_j T_j \mathcal{B}_j(f) = \sum_{j,k} \gamma_j(-u^1_k)u^2_k \o_{\rho(A)} \gamma_k(u^1_j f(u^2_j -)) = $$ 
$$ \sum_{j,k} \gamma_j(-)u^2_k \o \gamma_k(u^1_j f(u^2_j -) u^1_k = 
 \sum_j \gamma_j(-) u^1_j f(u^2_j) \o 1_{\ME} = f \o 1_{\ME}. $$

We also have a corollary for the right endomorphism algebra
$E := \End A_B$ extending $A$ via $\lambda: A \into E$.  

\begin{cor}
\label{cor-rightendo}
If $A \| B$ is a left D2 algebra extension, then $E \| A$
is a left D2, left balanced and a left Galois $S_{\rm cop}$-extension.
\end{cor}
\begin{proof}
We pass to the right D2 extension $A\op \| B\op$,
with left D2 extension, left balanced and left Galois
 $S_{\rm cop}$-extension $\rho: A\op \into \End {}_{B\op}A\op$
by the main theorem and Theorem~\ref{th-characterization}. Note again
 that $ \End {}_{B\op}A\op \cong E$ via $f \mapsto \overline{f}$ defined
in Prop.~\ref{prop-op} where $\overline{\lambda(a) f \lambda(a')}$ $=
\rho(\overline{a}) \overline{f} \rho(\overline{a'})$
for all $a,a'\in A$, $f \in E$, whence the isomorphic extension
$E \| A$ is left D2, left balanced and left Galois.
\end{proof}

The left $S$-coaction on $A$ implied by the corollary is given by
($g \in E$)
\begin{equation}
\label{eq: coact}
\delta_L(g) = \sum_i \beta_i \o_{R\op} g(-t^1_i)t^2_i
\end{equation}
which transferred to the smash product
$A \# S \cong \End A_B$ (via $a \# \alpha \mapsto \lambda(a) \circ \alpha$
\cite[3.8, 4.5]{KS}) comes out as 
\begin{equation}
\label{eq: smash prod}
\delta_L(a \# \alpha) = \alpha\2 \o a \# \alpha\1
\end{equation}
by comparing with eq.~(\ref{eq: left d2 comult}). \cite{Bo} briefly sketches how
a smash product $A \# S$ is a Galois extension over $A$ if $S$ (is a left
bialgebroid acting on $A$ and) possesses an antipode
(e.g., if $A \| B$ is H-separable).  

Every one-sided D2 extension has a canonical $T$- and $S$-bialgebroid associated to it;
the next corollary observes (thanks to B\"ohm) that the $S$-bialgebroid of a right D2 extension is isomorphic
to the $T\op$-bialgebroid of its left D2 endomorphism algebra extension.

\begin{cor}
Under the hypotheses in Theorem~\ref{th-endogalext},
the Galois mapping $\beta$ in eq.~(\ref{eq: gm}) restricts
to an isomorphism of left bialgebroids
\begin{equation}
\label{eq: tee-ess}
(\ME \o_{\rho(A)} \ME)^{\rho(A)} \cong S
\end{equation}
where the left-hand structure has the opposite multiplication of that in eq.~(\ref{eq: tee mult}). 
\end{cor}
\begin{proof}
The isomorphism of $K$-modules follows from noting that $\ME^{\rho(A)} = \End {}_BA_A \cong R$ via
$f \mapsto f(1)$, 
and applying the particular $\rho(A)$-$\rho(A)$-bimodule linearity of
the Galois mapping.  If $T^1 \o T^2 \in (\ME \o_{\rho(A)} \ME)^{\rho(A)} $,
then $\beta(T^1 \o T^2) = T^1(- T^2(1))$, which is in $\End {}_BA_B$ since
$T^1(-T^2(a)) = T^1(-T^2(1))a$ for each $a \in A$.  

The mapping $\beta$ is an algebra isomorphism since 
$$\beta(U^1 T^1 \o T^2 U^2) = U^1T^1(- T^2(U^2(1))) = $$
$$ U^1 (T^1(- T^2(1))U^2(1)) = \beta(U^1 \o U^2) \circ \beta(T^1 \o T^2),
$$
where also $U^1 \o U^2 \in (\ME \o_{\rho(A)} \ME)^{\rho(A)} $. 

The mapping $\beta$ commutes with the source, target, counit and comultiplication mappings of each bialgebroid;
e.g., $\beta( \id_A \o \lambda(r)) = \rho(r)$ for each $r \in R$,
and $\eps_S(\beta(T^1 \o T^2) = T^1(T^2(1))$ to which $\eps_T(T^1 \o T^2) = T^1 T^2$ maps
under the identification $\End {}_BA_A \cong R$. 
\end{proof}

There has been a question of whether a
 right or left progenerator H-separable extension $A \| B$ is split 
(i.e., has a $B$-$B$-bimodule projection $A \to B$), whence Frobenius: an affirmative answer implies some generalizations of results of Noether-Brauer-Artin on simple
algebras \cite{Su}. Unfortunately, the next example
of a one-sided free H-separable non-Frobenius extension
rules out this possibility.      
  
\begin{example}
\begin{rm}
Let $K$ be a field and $B$ the $3$-dimensional algebra of  upper triangular $2 \times 2$-matrices, which is not self-injective.
Since $B \| K1$ is trivially D2, the endomorphism algebra $A := \End B_K \cong M_3(K)$ is a left
D2 extension of $\lambda(B)$, which w.r.t.\ the ordered basis $\bra e_{11}, e_{12}, e_{22} \ket$
$\lambda(B)$ is the subalgebra of matrices
$$ [x,y,z] :=  \left( \begin{array}{ccc}
x & 0 & 0 \\
0 & x & y \\
0 & 0 & z
\end{array}
\right)
$$  The left D2 quasibases for $A \| B$ obtained from the opposite version of eq.~(\ref{eq: left D2 qb endo})
indicate an H-separable extension since the $\mathcal{B}_j$ are all of the form $\rho(r_j)$ for $r_j \in R$.
Now the centralizer $R$ is the $3$-dimensional algebra spanned by matrix units $e_{11}$, $e_{21}$, $e_{22} + e_{33}$.  
The module ${}_BA$ is free with left $B$-module isomorphism $A \to B^3$ given by ``separating out the columns''
$$(a_{ij}) \mapsto
([a_{11},a_{21},a_{31}], [a_{12},a_{22}, a_{32}], [a_{13}, a_{23}, a_{33}])$$
whence $A \o_B A$ and $\Hom (R_K,A_K)$ are both $27$-dimensional.  The $A$-$A$-homomorphism 
$A \o_B A \to \Hom (R_K, A_K)$ given by $a \o c \mapsto (r \mapsto arc)$ is easily computed to 
be surjective, therefore an isomorphism, whence ${}_A A \o_B A_A \cong {}_AA^3_A$, which shows
$A \| B$ is H-separable (and D2). By \cite[6.1]{M}, the extension $A \| B$ is not Frobenius
since $B$ is not a Frobenius algebra; therefore $A \| B$ is not split. (Alternatively, if
 there is a $B$-linear projection $E: A \to B$, we note $E(e_{32}) = 0$, so $e_{33} = e_{32} e_{23} \in \ker E$,
a contradiction.) 
 By applying the matrix transpose, the results of this example may be transposed to a right-sided version. 
\end{rm}
\end{example}

The endomorphism ring theorem for Galois extensions below 
follows directly from the main theorem and Theorem~\ref{th-characterization}.
Let $T'$ be  a right $R'$-bialgebroid with ${}_{R'}T'$ finite projective
and $S$ is as usual the left bialgebroid $\End {}_BA_B$ over the centralizer.

\begin{cor}[Endomorphism Ring Theorem for Galois Extensions]
If $A \| B$ is a right Galois $T'$-extension, then $\ME \| A\op$ is a left
$S$-Galois extension.  
\end{cor}

Next we note a converse of the endomorphism ring theorem in case
of Frobenius D2 extensions. Recall that for a Frobenius extension $A \| B$, being right D2  is equivalent
to being left D2 \cite{KS, BS}. Also recall our notation $E := \End A_B$.  We will essentially
note that under favorable conditions  equal modules  are cancellable in reversing the argument
in \cite[6.7]{KS}

\begin{theorem}
\label{th-converse-endo}
Let $A \| B$ be a right generator Frobenius extension.  If $E \| A$ is D2, then
$A \| B$ is D2.
\end{theorem}
\begin{proof}
Since $A \| B$ is Frobenius with Frobenius coordinate
system $\phi \in$ $ \Hom ({}_BA_B, {}_BB_B)$, $\sum_i x_i \o y_i \in (A \o_B A)^A$, then $A_B$ is finite projective, $E \cong A \o_B A$ via
$f \mapsto \sum_i f(x_i) \o y_i$ with inverse $a \o a' \mapsto$ $ \lambda(a) \circ \phi \circ \lambda(a')$, and $\lambda: A \into E$ is itself a Frobenius extension \cite{NEFE}.
Then by the hypothesis $A_B$ is a progenerator, so $E$ and $B$ are Morita equivalent with context bimodules
${}_EA_B$ and ${}_B\Hom (A_B, B_B)_E$, where we note that $\Hom (A_B, B_B) \o_E A \cong B$. 

If $E \| A$ is right D2, then $${}_EE \o_A E_A \oplus * \cong \oplus^n {}_EE_A,$$
where substitution of ${}_EE_A \cong {}_EA \o_B A_A$ and cancelling $A \o_A -$ yields
$$ {}_EA \o_B A \o_B A_A \oplus * \cong \oplus^n {}_E A \o_B A_A. $$
Tensoring this from the left by ${}_B\Hom (A_B, B_B) \o_E -$ and the Morita property yields
${}_BA \o_B A_A \oplus * \cong \oplus^n {}_B A_A,$ whence $A \| B$ is left D2.
\end{proof}

By \cite[6.1]{M}, it is enough  to assume in the theorem above that
 $E \| A$ is D2 Frobenius and $A_B$ a progenerator.  The theorem fully answers Question~1 at the end of the paper 
\cite{KN}:  the two conditions in the definition of depth two Frobenius extensions are
equivalent, not independent, since split extensions automatically satisfy
the generator condition (and cf.\ \cite[2.1(5)]{KK}).

Finally, we propose a natural problem from the point of view
of this paper and \cite{fer}:
is there an example of a left D2 extension which is not right D2?



\begin{thebibliography}{XXXXXX}
\begin{small}
\bibitem{BaS}{I. B\'{a}lint and K.~Szlach\'anyi,
Finitary Galois extensions over non-commutative bases,
KFKI preprint (2004), \texttt{RA/0412122}.}
\bibitem{BW}{T.~Brzezi\'nski and R.~Wisbauer,
\textit{Corings and Comodules}, LMS \textbf{309}, Cambridge University Press, 2003.}
\bibitem{Bo}{G.~B\"{o}hm,
Galois theory for Hopf algebroids, KFKI preprint (2004), \texttt{RA/0409513}.}
\bibitem{BS}{G.~B\"ohm and K.~Szlach\'anyi,
Hopf algebroids with bijective antipodes: axioms, integrals and duals,
\textit{J. Algebra} \textbf{274} (2004), 708--750.}
\bibitem{CM} A.~Connes and H.~Moscovici,
Differential cyclic cohomology and Hopf algebraic structures in transverse geometry, Essays on geometry and related topics, Vol.~1, 2, 217--255, Monogr.\ Enseign.\ Math., \textbf{38}, \textit{Enseignement Math.}, Geneva, 2001. 
\bibitem{DT}{Y. Doi and M.~Takeuchi,
Cleft comodule algebras for a bialgebra, \textit{Comm.\ Algebra} \textbf{14}
(1986), 801--818.}
\bibitem{NEFE}{L.~Kadison,
\textit{New Examples of Frobenius Extensions}, University Lecture Series \textbf{14}, AMS, Providence, 1999.
Update, 6 pp: www.ams.org/bookpages.}
\bibitem{fer}{L. Kadison,
A Galois theory for bialgebroids, depth two and normal Hopf subalgebras, in: Joint Proc.\ of the conf.\ Ferrara and Swansea
(2004), eds. T. Brezinski and C. Menini, Annali dell' Universita di Ferrara, Sez. V11,
Scienze Matematiche, to appear. \texttt{QA/0502188}.}
\bibitem{KK}{L. Kadison and B.~K\"ulshammer,
Depth two, normality and a trace ideal condition
for Frobenius extensions, \textit{Comm.\ Algebra}, to appear. \texttt{GR/0409346}.}
\bibitem{KN}{L.~Kadison and D.~Nikshych,
Hopf algebra actions on strongly separable extensions of depth two,
\textit{Adv.\ in Math.} \textbf{163} (2001), 258--286.}
\bibitem{KSt} L.~Kadison and A.A.~Stolin, 
             Separability and Hopf algebras,
             in: \textit{Algebra and its Applications},  
eds.\ D.V.~Huynh, S.K.~Jain, and S.~Lopez-Permouth, Contemporary  Math.\ \textbf{259}     
    A.M.S., Providence, 2000, 279-298. 
\bibitem{KS}{L.~Kadison and K.~Szlach\'anyi,
Bialgebroid actions on depth two extensions and duality, \textit{Adv.\
in Math.} \textbf{179} (2003), 75--121.}
\bibitem{KT}{L.~El Kaotit and J.~G\'omez-Torrecillas,
Comatrix corings: Galois corings, descent theory and
a structure theorem for cosemisimple corings,
\textit{Math.\ Z.} \textbf{244} (2003), 
887--906.}
\bibitem{Lu}{J.-H. Lu,
Hopf algebroids and quantum groupoids,
\textit{Int.\ J. Math.} \textbf{7} (1996), 47--70.}
\bibitem{Mo}{S.~Montgomery,
\textit{Hopf Algebras and Their Actions on Rings},  CBMS Regional Conf.\ 
Series in Math.\ 
Vol.\ 82, AMS, Providence, 1993.}
\bibitem{M}{K. Morita,
The endomorphism ring theorem for Frobenius extensions,
\textit{Math.\ Z.} \textbf{102} (1967), 385--404.}
\bibitem{Mu}{B.~M\"{u}ller, 
Quasi-Frobenius Erweiterungen I, {\it Math.\ Z.} {\bf 85} (1964), 345--368. 
Q.-F. E. II, ibid. {\bf 88} (1965), 380--409.} 
\bibitem{NV} D.~Nikshych and L.~Vainerman,
Finite dimensional quantum groupoids
and their applications, in:  \textit{New Directions in Hopf Algebras},
S.~Montgomery and H.-J.~Schneider, eds., MSRI Publications, vol.\ \textbf{43}, Cambridge, 2002, pp.\ 211--262. 
\bibitem{Su}{K.~Sugano,
private communication, 1999.}
\bibitem{Sz}{K.~Szlach\'anyi, Galois actions by finite quantum groupoids,
\textit{Locally compact quantum groups and groupoids (Strasbourg, 2002)}, IRMA Lect.\ Math.\ Phys.\ 2, de Gruyter, Berlin, 2003, 105--125. \texttt{QA/0205229}. }
\bibitem{T}{M.~Takeuchi, 
Groups of algebras over $A \o \overline{A}$, \textit{J. Math.\ Soc.\ Japan}
\textbf{29} (1977), 459--492.} 
\end{small}
\end{thebibliography}
\end{document}